\documentclass[11pt]{article}
\usepackage{amssymb, amsmath}
\newcommand{\s}{\,\,\,\,}
\newcommand{\bthm}[2]{\vskip 8pt\bf #1\hskip 2pt\bf#2\it \hskip 8pt}
\newcommand{\ethm}{\vskip 8pt\rm}

\def\dint{\displaystyle{\int}}
\def\lan{\langle}
\def\ran{\rangle}

\title{ A simplification of the proof of the existence of the extremal function for the
Moser-Trudinger inequality}
\author{Yuxiang Li
\\ {\it\small Math. Group, The abdus salam ICTP,Trieste 34100,
   Italy}
\\ {\it\small E-mail address:
liy@ictp.trieste.it}}
\date{}
\begin{document}
\maketitle
\section{Introduction:}
In \cite{L1} and \cite{L2}, the author has proved the existence of the extremal function
for
the Moser-Trudinger inequality on a compact Riemannian manifold. In \cite{L2},
one of the key proposition is the following

\bthm{Proposition}{1.1}  Let $(M,g)$ be a closed compact Riemannian manifold, and
$$F(u)=\dint_Me^{\alpha_n|u|^\frac{n}{n-1}}.$$
If the sup of the $F$ on
$$\mathcal{H}=\{u:\dint_M|\nabla u|^ndV_g=1,\s \dint_MudV_g=0\},$$
can not be attained, then
$$\sup_{u\in\mathcal{H}}F(u)
\leq |M| + \frac{\omega_{n-1}}{n} e^{\alpha_nS_p+1+1/2 + \cdots
+ 1/(n-1)}$$
for some $p\in M$.\ethm
Here $\alpha_n=n\omega_{n-1}^\frac{1}{n-1}$, $\omega_{n-1}$ is the volume
of the unit sphere on $\mathbb{R}^n$, and
$$S_p = \lim\limits_{x \rightarrow p}
(G+\frac{1}{\alpha_n} \log{dist(p,x)^n}).$$

In section 5 of  \cite{L2}, the author proved the above proposition
with capacity estimate. Since we had no idea on the speed of
convergence of the blow up sequence, we had to choose a sequence of
Green functions to get the estimate. So the proof seemed  too
complicated to be understood. In this short paper, we will simplify
the proof.  I hope the new proof
 is easier to be
understood .

The proof in this paper depends on the following theorem belong to
Carleson and Chang:

\bthm{Theorem}{A.1(\cite{C-C})} Let $B$ be the unit ball in $\mathbb{R}^n$.
Assume $u_k$ be a sequence in  $H^{1,n}_0(B)$ with $\int_B|\nabla u_k|^n=1$.
If $|\nabla u_k|^ndx\rightharpoondown \delta_0$, then
$$\limsup_{k\rightarrow+\infty}\dint_Be^{\alpha_n|u_k|^\frac{n}{n-1}}
\leq|B|(1+e^{1+1/2+\cdots+1/(n-1)}).$$\ethm

{\bf Remark:} Since we do not get any new results, we will not publish this paper.

\section{Review of results in section 4 of \cite{L2}}

Let $u_k\in H^{1,n}_0(M)$ which satisfies $\int_Mu_k=0$, $\int_M|\nabla u_k|^ndV_g=1$ and
$$\int_Me^{\beta_k|u_k|^\frac{n}{n-1}}
dV_g=\sup_{\int_M|\nabla v|^ndV_g=1,\int_Mv=0}
  \int_Me^{\beta_k|v|^\frac{n}{n-1}}dV_g,$$
where $\{\beta_k\}$ is an increasing sequence which converges to
$\alpha_n$. Then we have
$$-div|\nabla u_k|^{n-2}\nabla u_k=\frac{|u_k|^\frac{1}{n-1}sign(u_k)
}{\lambda_k}e^{\beta_k
  |u_k|^\frac{n}{n-1}}-\gamma_k,$$
where $sign(u_k)$ is the sign of $u_k$.
It is obvious that
$$\lambda_k=\dint_M|u_k|^{\frac{n}{n-1}}e^{\beta_k|u_k|^\frac{n}{n-1}}dV_g
 \hbox{,  }
\gamma_k=\dint_M\frac{|u_k|^\frac{1}{n-1}sign(u_k)}{\lambda_k}e^{
 \beta_k|u_k|^\frac{n}{n-1}}dV_g,$$
and
$$\sup_{\int_M|\nabla u|^ndV_g=1,\int_MudV_g=0}\dint_Me^{\alpha_n|u|^\frac{n}{
  n-1}}dV_g=\lim_{k\rightarrow\infty}\dint_Me^{\beta_k|u_k|^\frac{n}{
  n-1}}dV_g.$$

Let $c_k=\max\limits_{x\in M}|u_k(x)|=\max\limits_{x\in M}u_k(x)$, and $x_k\rightarrow p$.
Let $\{e_i(x)\}$ an  orthogonal
basis of $TM$ near $x_0$ and $exp_x: T_xM\rightarrow M$ be the exponential mapping.
 Let $B_r$ be the ball in $\mathbb{R}^n$. The
smooth mapping  $E: B_\delta(x_0)\times B_r\rightarrow M$ is define to be
$$E(x,y)=exp_x(y^ie_i(x)).$$

We set
$$g_{ij}(x,y)=\lan
(exp_x)_*\frac{\partial}{\partial y^i},(exp_x)_*\frac{\partial}{\partial y^j}
\ran_{E(x,y)}.$$
It is well-known that, we  can find a constant $a$, s.t.
$$\|g(x,y)-I\|_{C^0(B_\delta\times B_r)}\leq a|y|^2$$
when $\delta$ and $r$ are sufficiently small.

Let $r_k^n=\frac{\lambda_k}{c_k^\frac{n}{n-1}e^{\beta_kc_k^\frac{n}
{n-1}}}$, then we have $r_k\rightarrow 0$ and
$$\frac{n}{n-1}\beta_kc_k^\frac{1}{n-1}(u_k(E(x_k,r_kx))-c_k)\rightarrow
-n\log(1+c_nr^\frac{n}{n-1})$$
on any $B_L(0)$, where $c_n=(\frac{\omega_{n-1}}{n})^\frac{1}{n-1}$. Note that
$u_k(E(x_k,r_kx))$ is the function defined on $B_{\frac{r}{r_k}}\subset\mathbb{R}^n$.

Moreover, we can get $c_k^\frac{1}{n-1}u_k\rightharpoondown G$ in $H^{1,q}(M)$ for any
$q<n$, and
$$c_k^\frac{1}{n-1}u_k\rightarrow G\eqno (2.1)$$
$C^1$ smoothly on any $M\setminus B_\delta(p)$, where $G$
is the Green function  defined by
$$\left\{\begin{array}{l}
           -\Delta_nG=\delta_p-\frac{1}{\mu(M)}\\
           \int_MG=0.
          \end{array}\right.$$
We can prove that
$$\lim_{k\rightarrow +\infty}
\int_Me^{\beta_k|u_k|^\frac{n}{n-1}}
dV_g=\mu(M)+
\lim_{L\rightarrow+\infty}\lim_{k\rightarrow+\infty}
\int_{B_{Lr_k}(x_k)}e^{\beta_k|u_k|^\frac{n}{n-1}}
dV_g.\eqno (2.2)
$$
Moreover, on any domain $\Omega$ which contains $p$, we have
$$\lim_{L\rightarrow+\infty}\lim_{k\rightarrow+\infty}
\dint_{\Omega\setminus B_{Lr_k}(x_k)}e^{\alpha_ku_k^{\frac{n}{n-1}}}
=|\Omega|.\eqno (2.3)$$

\section{The proof of Proposition 1.1}
Let $\delta$ and $r$ be sufficiently small s.t.
$$|g_{ij}(x_k,y)-\delta_{ij}|\leq A|y|^2$$
for some constant $A$ which is independent of $k$.
Then, we can find a constant $a$, s.t.
$$|\nabla_0\varphi|^2=\sum_{i=1}^n|\frac{\partial u_i}
{\partial y^i}|^2\leq (1+a|y|^2)|\nabla\varphi|^2$$
for any $\varphi\in C^\infty$.

Let
$b_k=\sup_{\partial B_\delta}u_k$, and $u_k'=(u_k-b_k)^+$. We set
$f(y)=\frac{|g(x_k,y)|^\frac{1}{2n}}{\sqrt{1+a|y|^2}}$, where $a$ is chosen suitably large
so that $f\leq 1$ on $B_\delta(0)$.
We have
$$|\nabla_0fu_k'|\leq |\nabla u_k'||g|^\frac{1}{2n}+|\nabla_0 f||u_k'|.$$
Then
$$|\nabla_0fu_k'|^n\leq |\nabla u_k'|^n\sqrt{|g|}+c\sum_{k=1}^n|
\nabla u_k'|^{n-k}|u_k'|^{k}.$$
Here, $c$ is a constant.

Then
$$\dint_{B_\delta(x_k)}|\nabla_0 fu_k'|^ndx\leq 1-\dint_{M\setminus B_\delta(x_k)}|\nabla u_k'|^n
dV_g
+\frac{\rho(\delta)}{c_k^\frac{n}{n-1}},$$
where
$$\rho(\delta)=2c\sum_{k=1}^n\dint_{B_{2\delta}(p)}|\nabla G|^{k}|G|^{n-k}dV_g.$$

Given $\delta>0$, we set
$x_\delta\in \partial B_\delta$, s.t. $G(x_\delta)=\sup\limits_{x\in B_\delta}G(x)$.
Let
$$G(x_\delta)=-\frac{1}{\alpha_n}\log{\delta^n}+S_p+b_\delta,$$
and
$$\Omega_\delta=\{x:G(x)\geq G_\delta\}.$$
Then $b_\delta\rightarrow 0$ as $\delta\rightarrow 0$ and
$$\begin{array}{lll}
   \lim\limits_{k\rightarrow +\infty}
   \dint_{M\setminus B_\delta(x_k)}c_k^\frac{n}{n-1}|\nabla u_k'|^ndV_g&
            =&\dint_{M\setminus B_\delta(p)}|\nabla G|^ndV_g\\[1.7ex]
             &=&\dint_{\Omega_\delta}|\nabla G|^ndV_g\\[1.7ex]
            &=&\dint_{\Omega^c_\delta}\frac{G}{\mu(M)}dV_g+G_\delta
              \dint_{\partial\Omega_\delta}|\nabla G|^{n-2}
              \frac{\partial G}{\partial n}\\[1.7ex]
            &=&-2\dint_{\Omega_\delta}\frac{G}{\mu(M)}dV_g+G_\delta\\[1.7ex]
            &=&O(\delta^n\log{\delta})+G_\delta.
  \end{array}$$
If we set $\tau_k=\int_{B_\delta}|\nabla_0 fu_k'|^ndx$,  then by Theorem A.1., we have
$$\lim_{k\rightarrow+\infty}\dint_{B_\delta} e^{\frac{\beta_k}{\tau_k^{\frac{1}{n-1}}}
|fu_k'|^\frac{n}{n-1}}dx
\leq \delta^n\frac{\omega_{n-1}}{n}(1+e^{1+1/2+\cdots+1/(n-1)}).$$

On one hand, for any $\varrho<\delta$,
$$\lim_{L\rightarrow+\infty}\lim_{k\rightarrow+\infty}\dint_{B_{\varrho}
\setminus B_{Lr_k}(x_k)}e^{\frac{\beta_k}{\tau_k^\frac{1}{n-1}}|fu_k'|^\frac{n}{n-1}}
\leq O(\frac{1}{\delta^n})\lim_{L\rightarrow+\infty}\lim_{k\rightarrow+\infty}
\dint_{B_\varrho\setminus B_{Lr_k}(x_k)}
e^{\beta_k|u_k|^\frac{n}{n-1}}dV_g=o_\varrho(1).$$
On the other hand, for a fixed $\varrho$,
$$\lim_{k\rightarrow+\infty}\dint_{B_\delta\setminus\varrho}e^{\frac{\beta_k}
{\tau_k^\frac{1}{n-1}}|fu_k'|^\frac{n}{n-1}}dx\rightarrow |B_\delta\setminus B_\varrho|.$$
Letting $\varrho\rightarrow 0$, we get
$$\lim_{L\rightarrow+\infty}\lim_{k\rightarrow+\infty}\dint_{B_\delta
\setminus B_{Lr_k}(x_k)}e^{\frac{\beta_k}{\tau_k^\frac{1}{n-1}}|fu_k'|^\frac{n}{n-1}}=|B_\delta|.$$
Hence we get the following identity:
$$\lim_{L\rightarrow+\infty}\lim_{k\rightarrow+\infty}
\dint_{B_{Lr_k}(x_k)}e^{\frac{\beta_k}{\tau_k^\frac{1}{n-1}}|fu_k'|^\frac{n}{n-1}}\leq\delta^n
\frac{\omega_{n-1}}{n}e^{1+1/2+\cdots+1/(n-1)}.$$

If we fix an $L$, then  on $B_{Lr_k}(x_k)$, $f=1+O(e^{-\frac{\alpha_n}{2}c_k^\frac{n}{n-1}})$,
and then we have
$$\begin{array}{ll}
    \beta_k|u_k'|^\frac{n}{n-1}&=
       \frac{\tau_k^\frac{1}{n-1}}{f^\frac{n}{n-1}}\frac{\beta_k}{\tau_k^\frac{1}{n-1}}
         |fu_k'|^\frac{n}{n-1}\\
      &\leq \frac{\beta_k}{\tau_k^\frac{1}{n-1}}|fu_k'|^\frac{n}{n-1}
        (1-\frac{\frac{1}{n-1}G(x_\delta)-\epsilon(\delta,k)}{c_k^\frac{n}{n-1}})\\
      &=\frac{\beta_k}{\tau_k^\frac{1}{n-1}}|fu_k'|^\frac{n}{n-1}
        -\frac{\beta_k}{n-1}(\frac{u_k}{c_k})^\frac{n}{n-1}G(x_\delta)+\epsilon(\delta,k)
        |\frac{u_k}{c_k}|^\frac{n}{n-1}.
  \end{array}$$
where $\lim\limits_{\delta\rightarrow 0}
\lim\limits_{k\rightarrow+\infty}\epsilon(\delta,k)=0$.

Furthermore, we have
$$\beta_k|u_k|^\frac{n}{n-1}=\beta_k|u_k'+b_k|^\frac{n}{n-1}
\leq\beta_k|u_k'|^\frac{n}{n-1}+\beta_k\frac{n}{n-1}b_ku_k'+O(\frac{1}{c_k^\frac{1}{n-1}}).$$
It is easy to be deduced from (2.1) that
$$-\frac{\beta_k}{n-1}(\frac{u_k(x)}{c_k})^\frac{n}{n-1}G(x_\delta)
+\beta_k\frac{n}{n-1}b_ku_k'(x)\leq-\log{\delta^n}+\alpha_nS_p+o_k(1),$$
for any $x\in B_{Lr_k}(x_k)$ and sufficiently large $k$. Hence we get
$$\lim_{L\rightarrow+\infty}\lim_{k\rightarrow+\infty}
\dint_{B_{Lr_k}}e^{\beta_k|u_k|^\frac{n}{n-1}}\leq\frac{\omega_{n-1}}{n}
e^{1+1/2+\cdots+1/(n-1)+\alpha_nS_p}.$$
So, we can imply Proposition 1.1 from (2.2).

\end{document}